\documentclass[12pt]{article}
\usepackage{amsmath,amssymb,amsthm,enumerate}
\usepackage[cp1251]{inputenc}
\usepackage[english, russian, ukrainian]{babel}
\usepackage[T2A]{fontenc}
\usepackage[all]{xy}
\begin{document}
УДК 517.9
\bigskip

\textbf{Р. Р. Салимов, Б. А. Клищук} (R. R. Salimov, B. A.
Klishchuk), Институт  математики НАН Украины, Киев.
\bigskip

\textbf{НИЖНИЕ ОЦЕНКИ ДЛЯ ПЛОЩАДИ ОБРАЗА КРУГА}

(\textbf{LOWER BOUNDS FOR AREAS OF IMAGES OF DISCS}).

\bigskip

In this article we consider  $Q$-homeomorphisms  with respect to the
p-modulus on the complex plane with $p>2$. It is obtained a lower
area estimate for image of  discs  under such mappings. We solved
the extremal problem about minimization of the area functional of
images of discs.

\medskip

В работе рассматриваются $Q$-гомеоморфизмы относительно $p$-модуля
на комплексной плоскости при $p>2$. Получена нижняя оценка площади образа круга при таких
отображениях. Решена экстремальная проблема о минимизации функционала площади образа круга.

\newpage

{\bf 1. Введение.} Задача об искажении площадей при квазиконформных отображениях берет свое  начало в работе Б.Боярского, см. \cite{Bo}. Ряд результатов в этом направлении получен в работах \cite{GR}, \cite{Ast},  \cite{EH}.

Впервые верхняя оценка площади образа круга  при квазиконформных отображениях  встречается в монографии
М.А. Лаврентьева, см. \cite{Lav}. В монографии \cite{BGMR}, см. предложение 3.7,      получено  уточнение  неравенства   Лаврентьева в терминах угловой дилатации. Также ранее в работах  \cite{LS} и   \cite{S1} были получены верхние оценки искажения площади круга для  кольцевых и  нижних    $Q$-гомеоморфизмов. В данной работе получены нижние оценки площади образа
круга при $Q$-гомеоморфизмах относительно $p$-модуля при $p>2$.

Для простоты изложения ограничимся только плоским
случаем. Напомним некоторые определения. Пусть задано семейство $\Gamma$ кривых  $\gamma$ в
комплексной плоскости ${\Bbb C}$. Борелевскую функцию $\varrho:{\Bbb
C}\to[0,\infty]$ называют {\it допустимой} для $\Gamma$, пишут
$\varrho\in{\rm adm}\,\Gamma$, если
\begin{equation}\label{eq1.2KR} \int\limits_{\gamma}\varrho(z)\,|dz|\
\geqslant\ 1\qquad\forall\ \gamma\in\Gamma.\end{equation}
Пусть $p\in (1,\infty)$. Тогда   {\it
$p$--модулем} семейства $\Gamma$ называется величина
\begin{equation}\label{eq1.3KR}\mathcal{M}_p(\Gamma)\ =\
\inf_{\varrho\in\mathrm{adm}\,\Gamma}\int\limits_{{\Bbb
C}}\varrho^p(z)\,dm(z)\, .\end{equation}

Предположим, что  $D$ --- область в комплексной
плоскости $\mathbb{C}$, т.е. связное открытое подмножество
$\mathbb{C}$ и $Q: D\rightarrow [0,\infty]$
--- измеримая функция. Гомеоморфизм $f: D \rightarrow \mathbb{C}$
будем называть $Q$-гомеоморфизмом относительно $p$-модуля, если

\begin{equation}\label{eq1.4KR}\mathcal{M}_p(f\Gamma)\ \leqslant\
\int\limits_{{D}}Q(z)\,\varrho^{p}(z)\,dm(z)\, \end{equation} для
любого семейства $\Gamma$ кривых в $D$ и любой допустимой функции
$\varrho$ для $\Gamma$.

\medskip

Исследование  неравенств типа (\ref{eq1.4KR}) при $p=2$ восходит к
Л.~Альфорсу (см., напр.,  теорему 3, разд.~D, гл.~I,  \cite{A}), а
также О.~Лехто и К.~Вертанену (см. неравенство~(6.6), разд.~6.3, гл.~V в   \cite{LV}). В работе В.Я.~Гутлянского
(совместно с К.~Бишопом, О.~Мартио и М.~Вуориненом)  доказан многомерный  аналог  неравенства
(\ref{eq1.4KR})  для  квазиконформных
отображений  (см. \cite{BiGMV}).

\medskip

 Отметим также, что если
в (\ref{eq1.4KR}) функцию Q считать ограниченной п.в. некоторой постоянной $K\in  [1,\infty)$ и  $p = 2$, то мы приходим к классическим квазиконформным отображениям, которые были впервые введены в работах Грётча, Лаврентьева и Морри.


\bigskip

Пусть $Q:D\to[0,\infty]$ --- измеримая функция. Для  любого  числа $r>0$ обозначим

$$q_{z_0}(r)=
\frac{1}{2\pi\, r}\int\limits_{S(z_0, r)}Q(z)\, |dz| $$
--- среднее интегральное значение  функции $Q$ по  окружности $S(z_0,r)=\{z\in
\mathbb{C}: \,  |z-z_0| = r\}$.

\medskip

{\bf Теорема~1.} {\it Пусть $D$ и $D'$ ---
ограниченные области в $\Bbb C$  и  $f: D \rightarrow D' $ ---
$Q$-гомеоморфизм относительно $p$-модуля, $p>2$, $Q \in L^{1}_{\rm
loc}(D \setminus\{z_{0}\})$. Тогда при всех $r\in (0, d_0)$, $d_0 = {\rm
dist}(z_{0}, \partial D)$, имеет место оценка

\begin{equation}
\label{a1*}
 |fB(z_0,r)| \geqslant
\pi\,\left(\frac{p-2}{p-1}\right)^{\frac{2(p-1)}{p-2}}
\left(\int\limits_{0}^{r}\frac{dt}{t^\frac{1}{p-1}\,q_{z_0}^{\frac{1}{p-1}}(t)}\right)^{\frac{2(p-1)}{p-2}}\,,
\end{equation}
где $B(z_0,r)=\{z\in \mathbb{C}: |z-z_0|\leqslant r\}\,$.}
\medskip

Отметим, что при $p>2$ и $Q(z)\leqslant K$ из теоремы 1 мы приходим
к результату для круга из работы \cite{Ge}, см. лемму 7.


\medskip



\bigskip

{\bf 3. Доказательство основной теоремы 1.} Приведем некоторые
вспомогательные сведения о емкости конденсатора. Следуя работе
\cite{MRV}, пару $\mathcal{E}=(A,C)$, где $A\subset\mathbb{C}$ ---
открытое множество и $C$ --- непустое компактное множество,
содержащееся в $A$, называем {\it конденсатором}. Конденсатор
$\mathcal{E}$ называется  {\it кольцевым конденсатором}, если
$\mathfrak{R}=A\setminus C$ --- кольцевая область, т.е., если
$\mathfrak{R}$ --- область, дополнение которой
$\overline{\mathbb{C}}\setminus \mathfrak{R}$ состоит в точности из
двух компонент. Конденсатор $\mathcal{E}$ называется {\it
ограниченным конденсатором}, если множество $A$ является
ограниченным. Говорят также, что конденсатор $\mathcal{E}=(A,C)$
лежит в области $D$, если $A\subset D$. Очевидно, что если
$f:D\to\mathbb{C}$ --- непрерывное, открытое отображение и
$\mathcal{E}=(A,C)$ --- конденсатор в $D$, то $(fA,fC)$ также
конденсатор в $fD$. Далее $f\mathcal{E}=(fA,fC)$.

\medskip

Пусть $\mathcal{E}=(A,C)$ --- конденсатор. Обозначим  через $\mathcal{C}_0(A)$ множество
непрерывных функций $u:A\to\mathbb{R}^1$ с компактным носителем.
$\mathcal{W}_0(\mathcal{E})=\mathcal{W}_0(A,C)$ --- семейство неотрицательных функций
$u:A\to\mathbb{R}^1$ таких, что 1) $u\in \mathcal{C}_0(A)$, 2)
$u(x)\geqslant1$ для $x\in C$ и 3) $u$ принадлежит классу ${\rm
ACL}$.
При $p \geqslant1$ величину
\begin{equation}\label{eqks2.6}{\rm cap}_p\,\mathcal{E}={\rm cap}_p\,(A,C)=\inf\limits_{u\in \mathcal{W}_0(\mathcal{E})}\,
\int\limits_{A}\,\vert\nabla u\vert^p\,dm(z)\,,\end{equation}
где
\begin{equation}\label{eqks2.5}\vert\nabla
u\vert= \sqrt{\left(\frac{\partial u}{\partial
x}\right)^2+\left(\frac{\partial u}{\partial y}\right)^2}\,
\end{equation}
называют
{\it $p$-ёмкостью} конденсатора $\mathcal{E}$. В дальнейшем мы будем
использовать установленное в работе
\cite{Sh} равенство
\begin{equation}\label{EMC}
{\rm cap}_p\,\mathcal{E}=\mathcal{M}_p(\Delta(\partial A,\partial C; A\setminus
C)),\end{equation} где для множеств $\mathcal{F}_1$, $\mathcal{F}_2$
и $\mathcal{F}$ в $\mathbb{C}$,
$\Delta(\mathcal{F}_1,\mathcal{F}_2;\mathcal{F})$ обозначает
семейство всех непрерывных кривых, соединяющих $\mathcal{F}_1$ и
$\mathcal{F}_2$ в $\mathcal{F}$\,.

Известно, что при $p\geqslant1$,  см. предложение 5 из \cite{Kru},
\begin{equation}\label{eqks2.8} {\rm
cap}_p\,\mathcal{E}\geqslant\frac{\left[\inf
l(\sigma)\right]^p}{|A\setminus C|^{p-1}}\,.
\end{equation}
Здесь $l(\sigma)$ --- длина гладкой  (бесконечно дифференцируемой)  кривой $\sigma$, которая является  границей $\sigma=\partial U$ ограниченного открытого
множества $U$, содержащего $C$ и содержащегося вместе со своим
замыканием $\overline{U}$ в $A$, а точная нижняя грань берется по
всем таким $\sigma$.

\medskip

{\it Доказательство теоремы 1.}
Пусть $\mathcal{E}=\left(A, C\right)$ --- конденсатор, где $A=\{z\in
D: |z-z_{0}|<t+\Delta t\}$, $C=\{z\in D: |z-z_{0}|\leqslant t\}$,\
$t+\Delta t < d_{0}$. Тогда $f\mathcal{E} = \left(fA,fC\right)$ ---
кольцевой конденсатор в $D^{\prime}$ и согласно (\ref{EMC}) имеем
равенство

\begin{equation}
\label{a2} {\rm cap}_{p}\,  f\mathcal{E}= \mathcal{M}_{p}\left(\Delta(\partial fA, \partial fC; f(A\setminus C)\right).
\end{equation}

В силу неравенства (\ref{eqks2.8}) получим
\begin{equation}
\label{a3}{\rm cap}_{p}\,  f\mathcal{E} \geqslant \frac{\left[\inf\
l (\sigma)\right]^{p}}{|fA\setminus fC|^{p-1}}\,.
\end{equation}
Здесь $l(\sigma)$ --- длина гладкой  (бесконечно дифференцируемой)  кривой $\sigma$, которая является  границей $\sigma=\partial U$ ограниченного открытого
множества $U$, содержащего $C$ и содержащегося вместе со своим
замыканием $\overline{U}$ в $A$, а точная нижняя грань берется по
всем таким $\sigma$.

С другой стороны, в силу определения $Q$-гомеоморфизма относительно
$p$-модуля, имеем

\begin{equation}
\label{a4} {\rm cap}_{p}\,  f\mathcal{E} \leqslant \int\limits
_{D}Q(z)\,\varrho ^{p}(z)\,dm(z)
\end{equation}
для любой $\varrho\in{\rm adm}\ \Delta(\partial A, \partial C; A\setminus C).$

Легко проверить, что  функция $$\varrho(z)= \begin{cases} \frac{1}{|z-z_{0}|\,\ln\frac{t+\Delta t}{t}},&  z \in A\setminus C\\
0,&  z \not\in A\setminus C \end{cases}$$ является допустимой для
семейства $\Delta(\partial A, \partial C; A\setminus C)$ и поэтому

\begin{equation}
\label{a5} {\rm cap}_{p}\,  f\mathcal{E} \leqslant
\frac{1}{\ln^{p}\left(\frac{t+\Delta t}{t}\right)}
\int\limits_{R}\frac{Q(z)}{|z-z_{0}|^{p}}\,dm(z),
\end{equation}
где $R = \{z\in D: t\leqslant |z-z_{0}|\leqslant t+\Delta t\}$.

Комбинируя неравенства (\ref{a3}) и (\ref{a5}), получим

\begin{equation}\label{a6}
\frac{\left[\inf\ l (\sigma)\right]^{p}}{|fA\setminus fC|^{p-1}}
\leqslant \frac{1}{\ln^{p}\left(\frac{t+\Delta
t}{t}\right)}\int\limits _{R}\frac{Q(z)}{|z-z_{0}|^{p}}\,dm(z).
\end{equation}
По теореме Фубини имеем
\begin{equation}
\label{a7}  \int\limits _{R}\frac{Q(z)}{|z-z_{0}|^{p}}\,dm(z) =
\int\limits_{t}^{t+\Delta
t}\frac{d\tau}{\tau^{p}}\int\limits_{S(z_{0},\tau)}Q(z)\,|dz| = 2\pi
\int\limits_{t}^{t+\Delta t} \tau^{1-p}\,q_{z_{0}}(\tau)\,d\tau,
\end{equation}
где $q_{z_{0}}(\tau) = \frac{1}{2\pi \tau}\, \int\limits_{S(z_0,
\tau)}\, Q(z)\, |dz|$ и $S(z_0,\tau)=\{z\in \mathbb{C}: |z-z_0| =
\tau\}$. Таким образом,

\begin{equation}
\label{a8} \inf\ l (\sigma) \leqslant
(2\pi)^{\frac{1}{p}}\,\frac{|fA\setminus
fC|^{\frac{p-1}{p}}}{\ln\left(\frac{t+\Delta
t}{t}\right)}\left[\int\limits_{t}^{t+\Delta t}
\tau^{1-p}\,q_{z_{0}}(\tau)\,d\tau\right]^{\frac{1}{p}}.
\end{equation}

Далее, воcпользовавшись изопериметрическим неравенством
\begin{equation}
\label{a9} \inf\ l(\sigma) \geqslant 2 \sqrt{\pi
|fC|},
\end{equation}
получим
\begin{equation}
\label{a10} 2 \sqrt{\pi \, |fC|}
\leqslant(2\pi)^{\frac{1}{p}}\,\frac{|fA\setminus
fC|^{\frac{p-1}{p}}}{\ln\left(\frac{t+\Delta
t}{t}\right)}\left[\int\limits_{t}^{t+\Delta t}
\tau^{1-p}\,q_{z_{0}}(\tau)\,d\tau\right]^{\frac{1}{p}}.\end{equation}

Определим функцию $\Phi(t)$ для данного гомеоморфизма $f$ следующим
образом
\begin{equation}
\label{b1}
\Phi(t) = |fB(z_{0},t)|,
\end{equation}
где $B(z_{0},t) = \{z\in
\mathbb{C}: |z - z_{0}|\leqslant t\}$. Тогда из соотношения
(\ref{a10}) следует, что
\begin{equation}
\label{a11} 2 \sqrt{\pi \, \Phi(t)} \leqslant
(2\pi)^{\frac{1}{p}}\,\frac{[\frac{\Phi(t+\Delta t)-\Phi(t)}{\Delta
t}]^{\frac{p-1}{p}}}{\frac{\ln(t+\Delta t)-\ln t}{\Delta
t}}\left[\frac{1}{\Delta t}\int\limits_{t}^{t+\Delta t}
\tau^{1-p}\,q_{z_{0}}(\tau)\,d\tau\right]^{\frac{1}{p}}.
\end{equation}

Устремляя в неравенстве (\ref{a11}) $\Delta t\to 0$, и
учитывая монотонное возрастание функции $\Phi$ по $t \in (0,d_{0})$,
для п.в. $t$ имеем:

\begin{equation}
\label{a12}
\frac{2\pi^{\frac{p-2}{2(p-1)}}}{t^\frac{1}{p-1}\,q_{z_{0}}^{\frac{1}{p-1}}(t)}
\leqslant \frac{\Phi'(t)}{\Phi^{\frac{p}{2(p-1)}}(t)}.
\end{equation}
Отсюда легко вытекает следующее неравенство:
\begin{equation}
\label{a13}
\frac{2\pi^{\frac{p-2}{2(p-1)}}}{t^\frac{1}{p-1}\,q_{z_{0}}^{\frac{1}{p-1}}(t)}\,
\leqslant
\left(\frac{\Phi^{\frac{p-2}{2(p-1)}}(t)}{\frac{p-2}{2(p-1)}}\right)^{\prime}.
\end{equation}

Поскольку $p>2$, то функция
$g(t)=\frac{\Phi^{\frac{p-2}{2(p-1)}}(t)}{\frac{p-2}{2(p-1)}}$
является  неубывающей на $(0, d_0)$, где $d_0 ={\rm dist}(z_{0},
\partial D)$. Интегрируя обе части неравенства по  $t \in
[\varepsilon,r]$  и учитывая, что
\begin{equation}
\int\limits_{\varepsilon}^{r} \left(\frac{\Phi^{\frac{p-2}{2(p-1)}}(t)}{\frac{p-2}{2(p-1)}}\right)^{\prime} dt =\int\limits_{\varepsilon}^{r} \ g'(t)\, dt \leqslant g(r)-g(\varepsilon) \leqslant
\frac{\Phi^{\frac{p-2}{2(p-1)}}(r)-\Phi^{\frac{p-2}{2(p-1)}}(\varepsilon)}{\frac{p-2}{2(p-1)}}\,,
\end{equation}
см., напр., теорему IV. 7.4  в   \cite{Sa}, получаем
\begin{equation}\label{a14}
2\pi^{\frac{p-2}{2(p-1)}}\,\int\limits_{\varepsilon}^{r}\frac{dt}{t^\frac{1}{p-1}q_{z_{0}}^{\frac{1}{p-1}}(t)}
\leqslant
\frac{\Phi^{\frac{p-2}{2(p-1)}}(r)-\Phi^{\frac{p-2}{2(p-1)}}(\varepsilon)}{\frac{p-2}{2(p-1)}}\,.
\end{equation}
Устремляя в неравенстве (\ref{a14}) $\varepsilon\to 0$,  приходим к
оценке
\begin{equation}
\label{a15} \Phi(r) \geqslant \pi\,
\left(\frac{p-2}{p-1}\right)^{\frac{2(p-1)}{p-2}}
\left(\int\limits_{0}^{r}\frac{dt}{t^\frac{1}{p-1}\,q_{z_{0}}^{\frac{1}{p-1}}(t)}\right)^{\frac{2(p-1)}{p-2}}.
\end{equation}
Наконец, обозначая в последнем неравенстве $\Phi(r) =
|fB(z_{0},r)|$, имеем
\begin{equation}
\label{a16} |fB(z_{0},r)| \geqslant
\pi\,\left(\frac{p-2}{p-1}\right)^{\frac{2(p-1)}{p-2}}
\left(\int\limits_{0}^{r}\frac{dt}{t^\frac{1}{p-1}\,q_{z_{0}}^{\frac{1}{p-1}}(t)}\right)^{\frac{2(p-1)}{p-2}}
\end{equation}
и тем самым завершаем доказательство теоремы 1.

\bigskip

{\bf 3. Следствия из теоремы 1.} Из теоремы 1 непосредственно
вытекают следующие утверждения.
\medskip

Воспользовавшись условием  $q_{z_{0}}(t) \leqslant q_{0}\,t^{-\alpha}$, оценим правую часть
неравенства (\ref{a1*})   и проведя элементарные преобразования
приходим  к следующему результату.

\medskip

{\bf Следствие 1.} {\it Пусть $D$ и $D'$ ---
ограниченные области в $\Bbb C$  и  $f: D \rightarrow D' $ ---
$Q$-гомеоморфизм относительно $p$-модуля при $p>2$. Предположим, что функция $Q$ удовлетворяет условию

\begin{equation}
\label{b2J}
q_{z_{0}}(t) \leqslant q_{0}\,t^{-\alpha},\, q_{0} \in (0, \infty)\,,\, \alpha  \in [0, \infty)
\end{equation}
для $z_{0}\in D$ и  п.в. всех $t\in (0, d_0)$, $d_0 = {\rm
dist}(z_{0}, \partial D)$. Тогда при всех $r\in (0, d_0)$ имеет место оценка

\begin{equation}
\label{a1p}
 |fB(z_0,r)| \geqslant
\pi^{-\frac{\alpha}{p-2}}\, \left(\frac{p-2}{\alpha+p-2}\right)^{\frac{2(p-1)}{p-2}}q_{0}^{\frac{2}{2-p}} \, |B(z_0,r)|^{1+\frac{\alpha}{p-2}}\,.
\end{equation}
\medskip
}

\medskip
В частности, полагая здесь   $\alpha=0$, получаем следующее
заключение.

\medskip

{\bf Следствие~2.} {\it Пусть $D$ и $D'$ ---
ограниченные области в $\Bbb C$  и  $f: D \rightarrow D' $ ---
$Q$-гомеоморфизм относительно $p$-модуля, $p>2$ и $q_{z_0}(t) \leqslant q_{0}<\infty$ для п.в. $t \in (0, \, d_0)$, $d_0 =  {\rm
dist}(z_{0}, \partial D)$.  Тогда
имеет место оценка
\begin{equation}
\label{a1}
 |fB(z_0, r)| \geqslant q_{0}^{\frac{2}{2-p}}\,|B(z_0,r)|
\end{equation}
для всех $r \in (0,d_0)$\,.
}

\medskip {\bf Следствие~3.} {\it Пусть выполнены условия теоремы~1 и
 $Q(z) \leqslant K<\infty$  для п.в. $z\in  D$. Тогда имеет место оценка

\begin{equation}
\label{a1}
 |fB(z_0,r)| \geqslant K^{\frac{2}{2-p}}\,|B(z_0,r)|.
\end{equation}
для всех $r \in (0,d_0)$\,.
}

\medskip {\bf Замечание 1.} Следствие~3 является частным случаем результата Геринга для $E = B(z_0,r)$, см. лемму 7 в \cite{Ge}.

\medskip



{\bf Следствие ~4.} {\it Пусть $f: \mathbb{B} \rightarrow \mathbb{B}$ ---
$Q$-гомеоморфизм относительно $p$-модуля при $p>2$. Предположим, что функция $Q(z)$ удовлетворяет условию

\begin{equation}
\label{b2}
q(t) \leqslant \frac{q_{0}}{t\,\ln^{p-1}\frac{1}{t}},\, q_{0} \in (0, \infty)\,,
\end{equation}
при п.в. всех $t\in (0, 1)$,  где $q(t)=\frac{1}{2\pi t}\, \int\limits_{S_{t}}\, Q(z)\, |dz|$
--- среднее интегральное значение над окружностью $S_{t}=\{z\in
\mathbb{C}: \,  |z| = t\}$. Тогда при всех $r\in (0, 1)$ имеет место оценка

\begin{equation}
\label{a1}
 |fB_{r}| \geqslant
\pi\, \left(\frac{p-2}{p-1}\right)^{\frac{2(p-1)}{p-2}}\, q_{0}^{\frac{2}{2-p}} \, \left(r\ln\frac{e}{r}\right)^{\frac{2(p-1)}{p-2}},\,
\end{equation}
\medskip
где $B_{r}=\{z\in \mathbb{C}: |z|\leqslant r\}.\,$
}
\medskip

\medskip
{\bf 4. Экстремальные задачи для функционала площади.}  Пусть  $Q: \mathbb{B}\rightarrow [0,\infty]$
--- измеримая функция, удовлетворяющая условию
\begin{equation}
\label{b20}
q(t) \leqslant q_{0}\,,\, q_{0} \in (0, \infty)
\end{equation}
при  п.в.  $t\in (0, 1)$,  где $q(t)=\frac{1}{2\pi t}\, \int\limits_{S_{t}}\, Q(z)\, |dz|$
--- среднее интегральное значение над окружностью $S_{t}=\{z\in
\mathbb{C}: \,  |z| = t\}$.

Пусть $\mathcal{H}=\mathcal{H}(q_{0}, p,  \mathbb{B})$ --- множество
всех $Q$-гомеоморфизмов  $f: \mathbb{B} \rightarrow \mathbb{C}$
относительно $p$-модуля при  $p>2$  с условием (\ref{b20}).
Рассмотрим на классе $\mathcal{H}$ функционал площади

\begin{equation}
\label{b3}
\mathbf{S}_{r}(f)= |fB_{r}|\,.
\end{equation}

\medskip {\bf Теорема~2.} {\it Для всех $r\in[0, 1]$ справедливо равенство
\begin{equation}
\label{b3} \min_{f \in \mathcal{H}} \mathbf{S}_{r}(f)= \pi \,
q_{0}^{\frac{2}{2-p}} \,   r^{2}\,.
\end{equation} }

\medskip

\medskip

{\it Доказательство.} В силу следствия 1 немедленно вытекает оценка

\begin{equation}\label{b00*} \mathbf{S}_{r}(f)\geqslant  \, \pi \,  q_{0}^{\frac{2}{2-p}} \,   r^{2} \,.
\end{equation}

Укажем  гомеоморфизм  $f \in \mathcal{H}$  на котором реализуется
минимум функционала $\mathbf{S}_r(f)$\,.  Пусть  $f_0:\mathbb{B} \to
\mathbb{C}$, где

\begin{equation}
f_0(z)=q_{0}^{\frac{1}{2-p}} \, z \,
\end{equation}

Очевидно, что равенство в  (\ref{b00*}) достигается на отображении $f_0$\,. Осталось показать, что отображение, определенное таким образом,
явля\-ет\-ся $Q$-гомеоморфизмом относительно $p$-модуля с
$Q(z)=q_0$. Действительно,

\begin{equation}
l(z,f_0)=L(z,f_0)= q_{0}^{\frac{1}{2-p}}\,, \, \quad J(z,f_0)=q_{0}^{\frac{2}{2-p}}
\end{equation}
и
\begin{equation}
 K_{I,\,p}(z,f_0)=\frac{J(z,f_0)}{l^p(z,f_0)}=q_0 \,.
\end{equation}

По теореме 1.1  из работы \cite{SS}  отображение   $f_0$ является   $Q$-гомеоморфизмом относительно $p$-модуля с
$Q(z)= K_{I,\,p}(z,f_0)=q_0$.

\bigskip



\medskip






\medskip

 \medskip

Авторы: \textbf{Руслан Радикович Салимов, Богдан Анатольевич Клищук}
\medskip

\textbf{Институт математики НАН Украины, Киев }

\medskip

E-mail: \textbf{ruslan623@yandex.ru, bogdanklishchuk@mail.ru}
\end{document}